\documentclass[12pt, reqno]{amsart}
\usepackage [latin1]{inputenc}
\usepackage{amscd}
\usepackage{epsfig}
\usepackage{amssymb}
\usepackage{amsmath}
\usepackage{amsthm}
\usepackage[T1]{fontenc}
\usepackage{ae,aecompl}

\usepackage {amsfonts} 
\usepackage {latexsym} 
\usepackage {exscale} 
\usepackage{amsmath}

\newcommand{\R} {\ensuremath{\mathbb{R}}}

\newcommand{\C} {\ensuremath{\mathbb{C}}}
\newcommand{\Z} {\ensuremath{\mathbb{Z}}}

\newcommand{\OO}{\mathcal{O}}
\renewcommand{\o}[1]{\overline{#1}}

\newcommand{\dq}{\overline{\partial}}
\newcommand{\wt}[1]{\widetilde{#1}}
\DeclareMathOperator{\Reg}{Reg}
\DeclareMathOperator{\Sing}{Sing}

\DeclareMathOperator{\Jac}{Jac}

\newtheorem {satz} {Satz} [section]
\newtheorem {lem} [satz] {Lemma}
\newtheorem {cor} [satz] {Corollary}
\newtheorem {defn} [satz] {Definition}

\newtheorem {thm} [satz] {Theorem}

\DeclareMathOperator{\supp}{supp}

\renewcommand{\theta}{\vartheta}

\begin{document}

\title[Dolbeault Complex with weights] 
{The Dolbeault Complex with weights according to normal crossings}

\author{J. Ruppenthal}

\address{Department of Mathematics, University of Wuppertal, Gau{\ss}str. 20, 42119 Wuppertal, Germany.}
\email{ruppenthal@uni-wuppertal.de}

\date{July 23, 2008}

\subjclass[2000]{32A10, 32C35, 32C36}
\keywords{Cauchy-Riemann equations, integral operators, Dolbeault complex.}

\begin{abstract}
In the present paper, we define a Dolbeault complex with weights according to normal crossings,
which is a useful tool for studying the $\dq$-equation on singular complex spaces by resolution of singularities
(where normal 
crossings appear naturally). The major difficulty is to prove that this complex 
is locally exact.
We do that by constructing a local $\dq$-solution operator which 
involves only Cauchy's Integral Formula (in one complex variable)
and behaves well for $L^p$-forms with weights according to normal crossings.
\end{abstract}



\maketitle

\section{Introduction}

The motivation for the present paper is as follows: one strategy to study the 
$\dq$-equation
on singular complex spaces is to use Hironaka's resolution of singularities in order to pull-back
the $\dq$-equation to a regular setting, where it is treatable much easier. 
One can achieve that the exceptional set
of such a desingularization consists of normal crossings only. 
It is therefore important to
study a Dolbeault complex with weights according to normal crossings. 
We will make that precise
in the following. A simple version of such a weighted Dolbeault complex was used in this setting in \cite{FOV1} and \cite{Rp4}.

\vspace{2mm}
Let $Y$ be an analytic variety in $\C^n$ of pure dimension $d$, and $\pi: M\rightarrow Y$
a resolution of singularities. So, we may assume that $M$ is a complex manifold of dimension $d$,
$\pi$ is a proper analytic map which is a biholomorphism outside the exceptional set
$$X=\pi^{-1}(\Sing Y),$$
and $X$ consists only of normal crossings (see \cite{AHL,BiMi,Ha}). Let $Y^*=\Reg Y$ carry the metric induced
by the canonical embedding $\iota: Y\hookrightarrow \C^n$, and let $M$ be given an arbitrary metric.
We denote by $dV_Y$ the volume element on $Y^*$, and by $dV_M$ the volume element on $M$.
We can assume that $\pi$ preserves orientation. Let $Q\in X$. Then there is a neighborhood $U$ of $Q$
in $M$ with local coordinates $z_1, ..., z_d$ such that we can assume $Q=0\in U\subset\C^d$,
and
$$X\cap U=\{z\in U: z_1\cdots z_m=0\}$$
for a certain integer $m$, $1\leq m\leq d$.\\

Then we can assume by Lemma 2.1 in \cite{Rp4} 
that there is a holomorphic function
$J\in \OO(U)$, vanishing exactly on $X\cap U$, such that 
\begin{eqnarray}\label{eq:jac}
|J|^2 = \det \Jac_\R \pi.
\end{eqnarray}
We may write $J=z^w=z_1^{w_1}\cdots z_d^{w_d}$, where $w=(w_1, ..., w_d)\in\Z^d$ is a multi-index
with $w_1, ..., w_m\geq 1$ and $w_{m+1}=...=w_d=0$.

\vspace{1mm}
Now, let $1\leq p\leq \infty$, and $f\in L^p(G)$ where
$G=\pi(U)^*=\Reg \pi(U)\subset Y^*$.
Then it follows from \eqref{eq:jac} that
\begin{eqnarray*}
\int_{U\setminus X} |\pi^* f|^p |J|^2 dV_M = \int_G |f|^p dV_Y,
\end{eqnarray*}
and this yields that (in multi-index notation for $|z|^{2w/p}$):
\begin{eqnarray*}
f\in L^p(G) \ \ \Leftrightarrow\ \ |J|^{2/p} \pi^* f = |z|^{2w/p} \pi^* f\in L^p(U).
\end{eqnarray*}

This gives reason to the following construction:

\begin{defn}
Let $D\subset \C^n$ be an open set, $1\leq p\leq \infty$ and $s=(s_1, ..., s_n)\in \R^n$ a real multi-index.
Then, we define:
$$|z|^s L^p_{0,q}(D):=\{ f \mbox{ measurable on } D: |z|^{-s} f \in L^p_{0,q}(D)\}.$$
$|z|^s L^p_{0,q}(D)$ is a Banach space with the norm 
$$\|f\|_{|z|^s L^p_{0,q}(D)}:= \||z|^{-s} f\|_{L^p_{0,q}(D)}.$$
We use the multi-index notation
$|z|^{-s}=|z_1|^{-s_1}\cdots |z_n|^{-s_n}$.
\end{defn}

The main objective of the present paper is to study the $\dq$-equation on $|z|^s L^p_{0,q}(D)$.
But that does not make sense in general for the usual $\dq$-operator.
It is therefore adequate to introduce the following weighted operator
($\dq$ has to be understood in the sense of distributions throughout the paper):

\begin{defn}\label{defn:dqk}
Let $k=(k_1, ..., k_n)\in \Z^n$ be an integer-valued multi-index,
and let $f$ be a measurable $(0,q)$-form on $D\subset \C^n$ such that
$$z^{-k} f \in L^1_{(0,q),loc}(D)\ \ \mbox{ and } \ \ \dq\big(z^{-k}f\big) \in L^1_{(0,q+1),loc}(D).$$
Then, we set
$$\dq_k f := z^k \dq\big( z^{-k} f\big) \in |z|^k L^1_{(0,q+1),loc}(D).$$
\end{defn}

Note that $\dq_k f=0$ exactly if $\dq(z^{-k}f)=0$. It is clear that $\dq_k\circ\dq_k=0$.
We will now use the abstract Theorem of de Rham in order to establish a link between the $\dq_k$-equation
$\dq_k g=f$ in $|z|^s L^p_{0,*}(D)$ and certain cohomology groups on $D$. The right coherent analytic sheaves
to look at are the following:

\begin{defn}
For $1\leq j\leq n$, let $\mathcal{I}_j=(z_j)$ be the sheaf of ideals of $\{z_j=0\}$ in $\C^n$ .
If $k=(k_1, ..., k_n)\in \Z^n$ is an integer-valued multi-index, let
$$\mathcal{I}^k \OO = \mathcal{I}_1^{k_1}\cdots \mathcal{I}_n^{k_n} \OO$$
as a subsheaf of the sheaf of germs of meromorphic functions.
\end{defn}

Note that we could as well consider the usual $\dq$-operator on sections
of a holomorphic line bundle $L^k=L_1^{k_1}\otimes \cdots \otimes L_n^{k_n}$
such that $\mathcal{I}^k \OO \cong \OO(L^k)$, where $\OO(L^k)$ is the sheaf of germs of
holomorphic sections in $L^k$.
This point of view is equivalent and wouldn't influence the presentation much.

\vspace{1mm}
We need to choose the right operator $\dq_k$ for given values of $p$ and $s$.
$k=k(p,s)$ should be the maximal value such that $|z|^s L^p_{loc} \subset |z|^k L^1_{loc}$.
It will become clear that this is a good choice at several points throughout the paper.
So:

\begin{defn}\label{defn:k}
Let $1\leq p \leq \infty$ and $s$ be real numbers. Then we call
\begin{eqnarray*}
k(p,s) := \max\{m\in\Z: |z_1|^s L^p_{loc}(\C) \subset |z_1|^m L^1_{loc}(\C)\}
\end{eqnarray*}
the $\dq$-weight of $(p,s)$. 
For $s=(s_1, ..., s_n)\in \R^n$, let
$k(p,s)=(k_1, ..., k_n)\in \Z^n$ be given by
$k_j := k(p,s_j)$, or, equivalently,
\begin{eqnarray*}
k(p,s) := \max\{m\in\Z^n: |z|^s L^p_{loc}(\C^n) \subset |z|^m L^1_{loc}(\C^n)\}.
\end{eqnarray*}
Then, we define for $0\leq q\leq n$ the sheaves $|z|^s \mathcal{L}^p_{0,q}$ by:
\begin{eqnarray*}
|z|^s \mathcal{L}^p_{0,q} (U) :=\{ f\in |z|^s L^p_{(0,q),loc}(U): \dq_{k(p,s)} f \in |z|^s L^p_{(0,q+1),loc}(U)\}
\end{eqnarray*}
for open sets $U\subset\C^n$ (it is a presheaf wich is already a sheaf).
\end{defn}

We will see later (Lemma \ref{lem:k1}) how the $\dq$-weight can be computed explicitely.
Now we can state the main result of the present paper:

\begin{thm}\label{thm:main}
For $1\leq p\leq \infty$ and $s\in \R^n$, let $k(p,s)\in \Z^n$ be the $\dq$-weight according to Definition \ref{defn:k}.
Then:
\begin{eqnarray}\label{eq:complex}
0 \rightarrow \mathcal{I}^k\OO \hookrightarrow
 |z|^s \mathcal{L}^p_{0,0} \xrightarrow{\ \dq_k\ }
 |z|^s \mathcal{L}^p_{0,1} \xrightarrow{\ \dq_k\ }
\cdots \xrightarrow{\ \dq_k\ }
 |z|^s \mathcal{L}^p_{0,n} \rightarrow 0
\end{eqnarray}
is an exact (and fine) resolution of $\mathcal{I}^k\OO$.
\end{thm}

By the abstract Theorem of de Rham, this implies that
$$H^q(U,\mathcal{I}^k\OO) \cong \frac{\mbox{ker }
(\dq_k: |z|^s \mathcal{L}^p_{0,q}(U) \rightarrow |z|^s \mathcal{L}^p_{0,q+1}(U))}
{\mbox{Im }
(\dq_k: |z|^s \mathcal{L}^p_{0,q-1}(U) \rightarrow |z|^s \mathcal{L}^p_{0,q}(U))}$$
for open sets $U\subset \C^n$. Thus, we can study the equation
$\dq_k g =f$
on $U$ by investigating the groups $H^q(U,\mathcal{I}^k \OO)$.
Due to the local nature of Theorem \ref{thm:main},
it is easy to deduce similar statements on complex manifolds,
which will be a helpful tool for studying the $\dq$-equation on singular spaces
as indicated in the beginning. We will do that in a second paper \cite{Rp7}.
Let us point out what is new about this approach.

\vspace{1mm}
For a complex projective variety $Z\subset\C\mathbb{P}^n$,
the Cheeger-Goresky-MacPherson conjecture (see \cite{CGM}) states that
the $L^2$-deRham cohomology  $H^*_{(2)}(Z^*)$
of the regular part of the variety $Z^*:=\mbox{Reg } Z$
with respect to the (incomplete) restriction of the Fubini-Study metric
is naturally isomorphic to the intersection cohomology of middle perversity $IH^*(Z)$
(which in turn is isomorphic to the cohomology of a small resolution of singularities).
Ohsawa proved this conjecture under the extra assumption that the variety has only isolated singularities (see \cite{Oh}),
while it is still open for higher-dimensional singular sets.
The early interest in the conjecture of Cheeger, Goresky and MacPherson was motivated in large
parts by the hope that one could then use the natural isomorphism and a Hodge decomposition
for $H^k_{(2)}(Z^*)$ to put a pure Hodge structure on the intersection cohomology of $Z$ (cf. \cite{CGM}). That
was in fact done by Pardon and Stern in the case of isolated singularities (see \cite{PaSt2}).
Their work includes the computation of the $L^2$-Dolbeault cohomology groups $H^{p,q}_{(2)}(Z^*)$
in terms of cohomology groups of a resolution of singularities (see also \cite{PaSt1}).

\vspace{1mm}
Let us now direct our attention to the case of Stein varieties.
Though one would expect similar relations in this (local) situation,
no such representation of the $L^2$-Dolbeault cohomology is known.
The best results include quite rough lower and upper bounds on the dimension
of some of the groups (see e.g. \cite{DFV}, \cite{Fo}, \cite{FOV2} or \cite{Rp4}).
The origin of the present work is the attempt to compute the $L^2$-Dolbeault cohomology groups
in the spirit of the work of Cheeger-Goresky-MacPherson, Ohsawa, Pardon-Stern and others
in terms of certain cohomology groups on a resolution of singularities.
But, in the absence of compactness, most of their arguments do not carry over to the local situation
and one has to develop 
new tools, one of them the Dolbeault complex with weights presented here.

\vspace{1mm}
In view of the large difficulties in computing the $L^2$-cohomology, it seems reasonable
to gain a broader view and better understanding by also considering $L^p$-Dolbeault cohomology groups
for arbitrary $1\leq p\leq \infty$. In fact, by use of the Dolbeault complex with weights,
is is for example possible to compute the $L^p$-Dolbeault cohomology groups on an irreducible homogeneous variety
with an isolated singularity $Y$ for all $p$ such that $2d/p\notin \Z$ (where $d=\dim Y$, see \cite{Rp7}) and for $p=1$.
This does not solve the $L^2$-problem but gives a quite precise idea what to expect for the $L^2$-groups.
A crucial point is that for the $L^p$-theory, we cannot use the well-known $L^2$-Hilbert space methods.
It seems reasonable to make use of integral formulas. This article is a first step in that direction.
Another application of integral formulas to singular spaces can be found in \cite{RuZe}.

\vspace{2mm}
For the proof of the main Theorem \ref{thm:main}, we need to solve the $\dq$-equation locally with
weights according to normal crossings. It is adequate to do this by using the inhomogeneous
Cauchy Integral Formula in one complex variable and to integrate just over lines parallel to the cartesian coordinates.
Following this idea, we were able to construct an integral solution operator (see Theorem \ref{thm:lp-solution}) satisfying the regularity properties
needed in Theorem \ref{thm:main}.

\vspace{1mm}
The paper is organized as follows: In section \ref{sec:CIF}, we present a version of Cauchy's Integral Formula (with weights)
which will fit our needs, and prove the relevant estimates for this integral operator (Theorem \ref{thm:weighted}).
In section \ref{sec:homotopy}, we give a $\dq$-homotopy formula by integral operators
which involve only iterated application of Cauchy's Integral Formula (Theorem \ref{thm:homotopy1}).

\vspace{1mm}
Modifying this construction, we obtain in section \ref{sec:solution}
a $\dq_k$-solution operator with weights according to normal crossings which involves only integration over complex lines (Theorem \ref{thm:lp-solution}).
This integral operator has strong regularity properties that are sufficient to prove Theorem \ref{thm:main} easily,
but it is of interest on its own 
(not only in the context of Theorem \ref{thm:main}).
In the last section, we complete the proof of Theorem \ref{thm:main}
and introduce another related complex which is of similar structure and use (Theorem \ref{thm:main'}).

\vspace{1mm}
\section{Weighted $L^p$-Regularity of Cauchy's Integral Formula}\label{sec:CIF}

Let $D\subset\subset\C$ be a bounded domain in the complex plane
and $k\in\Z$. For a measurable function $f$ on $D$, we define
\begin{eqnarray*}
{\bf I}_k^D f (z):= \frac{z^k}{2\pi i} \int_D f(\zeta) \frac{d\zeta\wedge d\o{\zeta}}{\zeta^k (\zeta-z)},
\end{eqnarray*}
provided, the integral exists. Note that
$$\frac{\partial}{\partial\o{z}} {\bf I}_k^D f = f$$
in the sense of distributions if $f(\zeta)/\zeta^k$ is integrable in $\zeta$ over $D$.
In this section, we will develop weighted $L^p$-estimates for these operators ${\bf I}_k^D$.
For the choice of the right weight factors in the occurring Cauchy Integrals,
we have to use the $\dq$-weight (see Definition \ref{defn:k}).

\begin{thm}\label{thm:weighted}
Let $D\subset\subset \C$ be a bounded domain, $1\leq p\leq\infty$ and $s$ real numbers,
and let $k=k(p,s)$ be the $\dq$-weight of $(p,s)$ according to Definition \ref{defn:k}.
Then ${\bf I}_{k}^D$ is a bounded linear operator
$${\bf I}_{k}^D: |z|^s L^p(D) \rightarrow |z|^{s+1-\epsilon} L^p(D)\ \ \mbox{ for all } \epsilon>0.$$
\end{thm}

Before we prove Theorem \ref{thm:weighted}, we need to know a little more about the $\dq$-weight.
It can be computed explicitly by the following:

\begin{lem}\label{lem:k1}
Let $1\leq p \leq \infty$ and $s$ be real numbers, and $k(p,s)$ the $\dq$-weight of $(p,s)$
according to Definition \ref{defn:k}. Then
\begin{eqnarray}\label{eq:k1}
k(p,s) = \left\{
\begin{array}{ll}
\max\{m\in\Z: m<2 + s-2/p\} & ,\ p\neq 1,\\
\max\{m\in\Z: m\leq 2 + s-2/p\}& ,\ p=1.
\end{array}\right.
\end{eqnarray}
\end{lem}

\begin{proof}
Let $k(p,s)$ be the $\dq$-weight of $(p,s)$, and let $k'(p,s)$ denote the right hand side in \eqref{eq:k1}.
By use of the H\"older Inequality we will now show firstly that \linebreak $k(p,s)\geq k'(p,s)$,
namely $|z|^s L^p_{loc}(\C) \subset |z|^{k'} L^1_{loc}(\C)$. So, let $f\in |z|^s L^p_{loc}(\C)$, $p>1$. Then,
for a bounded domain $D\subset\subset\C$:
\begin{eqnarray*}
\|z^{-k'} f\|_{L^1(D)} &\leq& \| |z|^{-s} f\|_{L^p(D)} \||z|^{s-k'} \|_{L^{\frac{p}{p-1}}(D)},
\end{eqnarray*}
which is bounded if $(s-k')\frac{p}{p-1} > -2$,
which is equivalent to $k'<2+s-2/p$.
If $p=1$, it is seen easily that we need $k'\leq s$.\\

It remains to show that $k(p,s)$ cannot be bigger than $k'(p,s)$.
So, we need to find functions in $|z|^s L^p_{loc}(\C)$ that are not in $|z|^{k'+1} L^1_{loc}(\C)$.
Note that
\begin{eqnarray*}
k'(p,s) = \left\{
\begin{array}{ll}
\min\{m\in\Z: m\geq1 + s-2/p\} & ,\ p\neq 1,\\
\min\{m\in\Z: m> 1 + s-2/p\}& ,\ p=1.
\end{array}\right.
\end{eqnarray*}
We will now distinguish two cases.
Let us first treat the situation that 
\begin{eqnarray}\label{eq:ineq1}
k'(p,s)>1+s-2/p.
\end{eqnarray}
Let $l:=s-2/p + \epsilon$. Then
$$|z|^l \in |z|^s L^p_{loc}(\C)$$
if $\epsilon>0$. On the other hand
$$|z|^{l} \notin |z|^{k'+1} L^1_{loc}(\C)$$
if $l-k'-1\leq-2$. But
$$(s-2/p) -k'-1 < -2$$
by assumption \eqref{eq:ineq1}. So, there exists $\epsilon>0$ such that $l-k'-1<-2$, implying
that $|z|^l \in |z|^s L^p_{loc}(\C)$ but $|z|^l \notin |z|^{k'+1} L^1_{loc}(\C)$.
Secondly, assume that 
$$k'(p,s)=1+s-2/p.$$
This implies that $p>1$. Consider
$$g(z) := |z|^{s-2/p} (\log |z|)^{-1}.$$
Then 
$$g\notin |z|^{k'+1} L^1_{loc}(\C),$$
because $|z|^{-2}(\log|z|)^{-1} \notin L^1_{loc}(\C)$, because
$(r \log r)^{-1} \notin L^1([0,1])$, having $\log\log r$ as antiderivate.
On the other hand,
$$g\in |z|^s L^p_{loc}(\C),$$
because
$|z|^{-2}\log^{-p} |z| \in L^1_{loc}(\C)$ for $p>1$, because $r^{-1} \log^{-p} r\in L^1([0,1])$,
having
$$(1-p)^{-1} \log^{1-p} r$$
as antiderivate for $p>1$.
\end{proof}

It is now easy to deduce furthermore:

\begin{lem}\label{lem:k2}
Let $1\leq p \leq \infty$ and $s$ be real numbers, and $k(p,s)$ the $\dq$-weight of $(p,s)$
according to Definition \ref{defn:k}. If $p<\infty$, then
\begin{eqnarray*}
k(p,s) &=& k_0(s)+k_1(p,s) = [s] + k_1(p,s)\\
 &:=& \max\{l\in\Z: l\leq s\} +\left\{
\begin{array}{ll}
2 &, \mbox{ if }\ (s-k_0)p>2,\\
1 &, \mbox{ if }\ 2 \geq (s-k_0)p > 2 - p,\\
0 &, \mbox{ if }\ 2-p \geq (s-k_0)p.
\end{array}\right.
\end{eqnarray*}
In case $p=\infty$, then 
$$k_1(\infty,s)=\left\{
\begin{array}{ll}
1 &, \mbox{ if } s-k_0=0,\\
2 &, \mbox{ if } s-k_0>0. 
\end{array}\right.$$
\end{lem}


\begin{proof}
The case $p=\infty$ follows directly from Lemma \ref{lem:k1}, because
$k(\infty,s) = [s] + 1$
if $s\in \Z$, and
$k(\infty,s)=[s]+2$
if $s\notin\Z$. If $p=1$, then 
$$k(1,s)=[s]=[s]+k_1(1,s),$$
because $k_1(1,s)=0$ for arbitrary $s$. If $p\notin\{1,\infty\}$, then
\begin{eqnarray*}
k(p,s) - [s] &=& \max\{l\in\Z: l+[s] < 2 +s -2/p\}\\
&=& \max\{l\in\Z: l < 2 +(s-[s]) -2/p\}\\
&=& \left\{
\begin{array}{ll}
2 &, \mbox{ if }\ 2< 2+(s-[s])-2/p,\\
1 &, \mbox{ if }\ 1 < 2+(s-[s])-2/p \leq 2,\\
0 &, \mbox{ if }\ 2+(s-[s])-2/p \leq 1.
\end{array}\right.\\
&=& \left\{
\begin{array}{ll}
2 &, \mbox{ if }\ (s-k_0)p>2,\\
1 &, \mbox{ if }\ 2 \geq (s-k_0)p > 2 - p,\\
0 &, \mbox{ if }\ 2-p \geq (s-k_0)p.
\end{array}\right.
\end{eqnarray*}
\end{proof}

There is another interesting concept in the context of our $\dq$-weight,
which we will discuss briefly before proving Theorem \ref{thm:weighted}:

\begin{defn}\label{defn:wtk}
Let $1\leq p \leq \infty$ and $s$ be real numbers. Then we call
$$\wt{k}(p,s):= \min\{m\in\Z: \mathcal{I}^m \OO_\C \subset |z|^s L^p_{loc}(\C)\}$$
the modified $\dq$-weight of $(p,s)$. It follows that
\begin{eqnarray}\label{eq:wtk}
\wt{k}(p,s) = \left\{
\begin{array}{ll}
\min\{k\in \Z: (s-k)p < 2\} &,\ p<\infty,\\
\min\{k\in \Z: (s-k) \leq 0\}&,\ p=\infty.
\end{array}\right.
\end{eqnarray}
\end{defn}

\begin{proof}
Let $k'=k'(p,s)$ denote the right hand side in \eqref{eq:wtk}.
By use of the H\"older Inequality, it is easy to see that $\wt{k}(p,s)\leq k'(p,s)$,
but it cannot be smaller, because 
$$z^{k'-1} \notin |z|^s L^p_{loc}(\C).$$
That can be seen as follows: If $p=\infty$, then \eqref{eq:wtk}
implies that $k'-1-s<0$.
If $p< \infty$, then \eqref{eq:wtk} implies that $(k'-1-s)p \leq -2$.
\end{proof}

The connection between $\dq$-weight and modified $\dq$-weight is clearyfied by:

\begin{lem}\label{lem:kk}
Let $1\leq p \leq \infty$ and $s$ be real numbers,
$k(p,s)$ the $\dq$-weight of $(p,s)$, and $\wt{k}(p,s)$ the modified $\dq$-weight.
Then
$$k(p,s) = \left\{
\begin{array}{ll}
[s]=s &, \mbox{ if } p=1 \mbox{ and } s\in\Z,\\
\min\{k\in\Z: (s-(k-1)) \leq 0 &,\ p=\infty,\\
\min\{k\in\Z: (s-(k-1))p\leq 2\} &, \mbox{ otherwise}.
\end{array}
\right.$$
It follows that $k(p,s)=\wt{k}(p,s)$
if $(s-[s])p\in\{2,2-p\}$, and
$k(p,s)=\wt{k}(p,s) + 1$
in all other cases.
\end{lem}

Note that this implies that there is an injective embedding $\mathcal{I}^k\OO \hookrightarrow |z|^s\mathcal{L}^p_{0,0}$
in the situation of Theorem \ref{thm:main}, because $k(p,s)\geq \wt{k}(p,s)$.

\begin{proof}
We will use the representation of $k(p,s)$ given in Lemma \ref{lem:k1} and Lemma \ref{lem:k2}.
The first two cases are immediate. So, assume that $p<\infty$,
and that $p>1$ or $s\notin\Z$. It is clear that
$$k':=k_0(s)+2 = [s] +2$$
satisfies $(s-(k'-1)) p \leq 2$
and $(s-(k'-1-3))p> 2$. Thus there exists $r(p,s)\in\{0,1,2\}$ such that:
$$\min\{k\in \Z: (s-(k-1))p \leq 2\} = k' - r(p,s).$$
It is not hard to see that $r(p,s)=2-k_1(p,s)$, and that proves the statement.
\end{proof}

We presented the modified $\dq$-weight at this place,
because it seems interesting to mention the following direct consequence of Theorem \ref{thm:weighted}:

\begin{cor}\label{cor:weighted}
Let $D\subset\subset \C$ be a bounded domain, $1\leq p\leq\infty$ and $s$ real numbers,
and let $\wt{k}=\wt{k}(p,s)$ be the modified $\dq$-weight of $(p,s)$ according to Definition \ref{defn:wtk}.
Then ${\bf I}_{\wt{k}}^D$ is a bounded linear operator
$${\bf I}_{\wt{k}}^D: |z|^s L^p(D) \rightarrow |z|^{s} L^p(D).$$
\end{cor}

\begin{proof}
The statement is clear if $\wt{k}(p,s)=k(p,s)$. If that is not the case,
we know that $(s-[s])p\notin\{2,2-p\}$ by Lemma \ref{lem:kk}.
But then it follows from the Definition of the $\dq$-weights and Lemma \ref{lem:kk}
that there exists $\epsilon>0$ such that $t:=s-1+\epsilon$ satisfies:
$k(p,t) = k(p,s) - 1 = \wt{k}(p,s)$.
And then it follows from Theorem \ref{thm:weighted} that the following map is bounded:
$${\bf I}_{\wt{k}}^D = {\bf I}^D_{k(p,t)}: |z|^s L^p(D) \hookrightarrow |z|^{t}L^p(D) \rightarrow |z|^{t+1-\epsilon} L^p(D) = |z|^s L^p(D).$$
\end{proof}

The rest of this section is dedicated to the proof of Theorem \ref{thm:weighted}.
It is based on two basic estimates of the integrals involved,
which we will recall for the convenience of the reader. Firstly, we need
(for a proof, see \cite{Rp2}, Lemma 6.1.1):

\begin{thm}\label{thm:cif1}
For $R>0$ and $0\leq\alpha,\beta<2$, let $\Delta_R:=\{z\in\C:|z|<R\}$ and
$$J_R(z):=\int_{\Delta_R}\frac{dV_\C(t)}{|t|^\alpha |t-z|^\beta}.$$
Then there exists a constant $C(\alpha,\beta)>0$,
depending only on $\alpha$ and $\beta$, such that
$$J_R(z)\leq C(\alpha,\beta) \left\{\begin{array}{llr}R^{2-\alpha-\beta} &, &\alpha+\beta<2,\\
1+\left|\log R-\log|z|\right| &, \mbox{ if } &\alpha+\beta=2,\\
|z|^{2-\alpha-\beta} &, &\alpha+\beta>2,\end{array}\right.$$
for all $z\in\C$ (with $z\neq 0$ if $\alpha+\beta\geq 2$).
\end{thm}

Moreover, we will make use of the following generalization of
the classical Young's inequality for convolution integrals (see e.g. \cite{Ra}, Appendix B, for the proof),
which is used frequently for estimating integral operators:

\begin{thm}\label{thm:young}
Let $(X,\mu)$ and $(Y,\nu)$ be two measure spaces, and suppose that $K$ is a measurable function on $X\times Y$
(with respect to product measure), which satisfies
$$\int_X |K(x,y)|^t d\mu(x)\leq M^t \ \ \mbox{ for almost all } y\in Y$$
and
$$\int_Y |K(x,y)|^t d\nu(y)\leq M^t \ \ \mbox{ for almost all } x\in X$$
for some $M<\infty$ and $s\geq 1$. Then the linear operator $f\mapsto {\bf T}f$ defined $\nu$-a.e. by
$${\bf T} f(y) = \int_X K(x,y) f(x) d\mu(x)$$
is bounded from $L^p(X)$ to $L^r(Y)$ with norm $\leq M$ for all $1\leq p,r\leq\infty$ with
$$\frac{1}{r}=\frac{1}{p} + \frac{1}{t} -1.$$
\end{thm}

Let use now show how we can use Theorem \ref{thm:cif1} and Theorem \ref{thm:young} to deduce
Theorem \ref{thm:weighted}. 
We will use the representation for the $\dq$-weight $k(p,s)$ given in Lemma \ref{lem:k2}.
It is easy to see that it is enough to consider the situation where $k_0(s)=[s]=0$ and $0\leq s<1$.
We have to distinguish four cases.
Firstly, let $k_1(p,s)=0$, hence
\begin{eqnarray}\label{eq:case1}
k(p,s)=0 \ \ \ \mbox{ and } \ \ \ (1+s)p\leq 2.
\end{eqnarray}

Let $1<t(\epsilon)<2$ be chosen quite close to $2$. This statement will be made precise later.
Then, let $r$ be defined by
\begin{eqnarray}\label{eq:hoelder1}
1/r = 1/p + 1/t -1.
\end{eqnarray}
Note that $r>p$ and $r\geq t$. We will now use the H\"older-Inequality with the three coefficients
$1/r$, $1/t-1/r$ and $1/p - 1/r$ (the sum equals $1$ by the use of \eqref{eq:hoelder1}):
\begin{eqnarray*}
\int_D |f(\zeta)| \frac{dV(\zeta)}{|\zeta-z|} &=&
\int_D \left(
\frac{|f(\zeta)|^p |\zeta|^{rs}}{|\zeta|^{sp}|\zeta-z|^t}\right)^{\frac{1}{r}}
\left(\frac{1}{|\zeta-z|^t}\right)^{\frac{1}{t}-\frac{1}{r}} \left(\left|\frac{f(\zeta)}{\zeta^s}\right|^p\right)^{\frac{1}{p}-\frac{1}{r}} dV(\zeta)\\
&\leq& C(t)^{1/t-1/r}\ \|f\|^{1-p/r}_{|z|^s L^p(D)} 
\left(\int_D \left|\frac{f(\zeta)}{\zeta^s}\right|^p \frac{|\zeta|^{rs} dV(\zeta)}{|\zeta-z|^t}\right)^{\frac{1}{r}},
\end{eqnarray*} 
where the second factor has been treated by Theorem \ref{thm:cif1} in order to receive the constant $C(t)>0$.
Here we have used the fact that $t$ is chosen to be $<2$.

Using this inequality, we can calculate (by the use of Fubini's Theorem and Theorem \ref{thm:cif1}):
\begin{eqnarray*}
&& \int_D \left||z|^{-s} {\bf I}_0^D f(z)\right|^r dV(z) \\ &\leq& 
C(t)^{r/t-1} \|f\|_{|z|^s L^p(D)}^{r-p} \int_D \int_D \left|\frac{f(\zeta)}{\zeta^s}\right|^p \frac{|\zeta|^{rs} dV(\zeta)}{|z|^{rs}|\zeta-z|^t} dV(z)\\
&=& C(t)^{r/t-1} \|f\|_{|z|^s L^p(D)}^{r-p} \int_D \left|\frac{f(\zeta)}{\zeta^s}\right|^p |\zeta|^{rs}
\int_D \frac{dV(z)}{|z|^{rs}|\zeta-z|^t} dV(\zeta)\\
&\leq& C(rs,t) C(t)^{r/t-1} \|f\|_{|z|^s L^p(D)}^{r-p} \int_D \left|\frac{f(\zeta)}{\zeta^s}\right|^p |\zeta|^{rs} |\zeta|^{2-rs-t} dV(\zeta)\\
&\lesssim& C(rs,t) C(t)^{r/t-1} \|f\|_{|z|^s L^p(D)}^{r-p} \|f\|_{|z|^s L^p(D)}^p \lesssim \|f\|_{|z|^s L^p(D)}^r.
\end{eqnarray*}
In order to apply Theorem \ref{thm:cif1}, we have used the fact that $rs<2$. 
That can be seen as follows:
The statement is trivial if $s=0$. If $s>0$, 
\eqref{eq:case1} implies that $p<2$ and that
$$sp \leq 2-p \ \ \Rightarrow \ \ sp/(2-p) \leq 1,$$
and from \eqref{eq:hoelder1} and $t<2$, we deduce that
$$1/r = 1/p+1/t -1 > 1/p -1/2 = (2-p)/2p.$$
Together, we receive:
$$rs < 2sp/(2-p) \leq 2.$$
Now we can complete the proof of the first case. Let $a(\epsilon)>2/p$ be chosen so small that $(1-\epsilon)ap<2$,
and let $b(\epsilon)<\frac{2/p}{2/p-1}=2/(2-p)$ be the dual coefficient such that $1/a+1/b=1$.
It is not hard to see that we can chose $t(\epsilon)$ such that $r=bp$, namely
$$t:= \frac{ap}{ap-1} < 2,$$
for this implies
$$\frac{1}{r}=\frac{1}{p} + \left(\frac{ap-1}{ap}\right) -1 
= \frac{1}{p} -\frac{1}{ap} = \frac{1}{p}\left(\frac{a-1}{a}\right) = \frac{1}{pb}.$$
So, by use of the H\"older-Inequality, $(1-\epsilon)ap<2$ and $|z^{-s} {\bf I}_0^D f|^{pb} = |z^{-s} {\bf I}_0^D f|^r$
(which has already been estimated), we finally get:
\begin{eqnarray*}
\|{\bf I}_0^D f\|^p_{|z|^{s+1-\epsilon} L^p(D)} &=&
\int_D \left||z|^{-1+\epsilon -s} {\bf I}_0^D f(z)\right|^p dV(z)\label{eq:ref1}\\
&\leq& \left(\int_D |z|^{(\epsilon-1)ap} dV(z)\right)^{1/a} \left( \int_D \left||z|^{-s} {\bf I}_0^D f(z)\right|^r dV(z)\right)^{1/b}\\
&\lesssim& \|f\|_{|z|^s L^p(D)}^{r/b} = \|f\|_{|z|^s L^p(D)}^p.\label{eq:ref3}
\end{eqnarray*}

For the case $k(p,s)=k_1(p,s)=1$, we have to distinguish two different situations.
Let us firstly consider the situation where $s>0$ and 
$$ 2 \geq sp > 2 - p.$$ 
This implies that $p\leq 2/s <\infty$.
It is enough to combine the principle of the last step in the case $k=k_1=0$ with Theorem \ref{thm:young}.
Actually, the proof of Theorem \ref{thm:young} is quite similar to the procedure that we have used
to estimate $\||z|^{-s}{\bf I}_0^D f\|_{L^r(D)}^r$. 
We can assume that $\epsilon<s$.
Choose 
$$a(\epsilon)>2/(ps)\geq 1$$
so small that $(s-\epsilon)ap<2$, and let $b< \frac{2/(ps)}{2/(ps)-1}=2/(2-ps)$ be the dual exponent such that
$1/a+1/b=1$. Then let $r=bp$ and choose $t$ such that
$$\frac{1}{r}=\frac{1}{p} + \frac{1}{t} -1.$$
Note that $r=bp>p$ implies that $t>1$. On the other hand, we see that 
$$\frac{1}{t} = \frac{1}{p}\left(\frac{1}{b}-1\right) + 1 > \frac{1}{p}\left(\frac{2-ps}{2} - 1\right) +1 = 1-\frac{s}{2},$$
which is equivalent to
\begin{eqnarray}\label{eq:ref4}
t (2-s) <2.
\end{eqnarray}

By use of the H\"older-Inequality, we deduce:
\begin{eqnarray*}
\|{\bf I}_1^D f\|^p_{|z|^{s+1-\epsilon} L^p(D)} &=&
\int_D \left||z|^{\epsilon-s-1} {\bf I}_1^D f\right|^p dV(z)\\
&\leq& \left(\int_D|z|^{(\epsilon-s)ap} dV(z)\right)^{1/a}
\left(\int_D \left||z|^{-1} {\bf I}_1^D f\right|^{bp} dV(z)\right)^{1/b}\\
&\lesssim& \||z|^{-1} {\bf I}_1^D f\|^p_{L^r(D)}.
\end{eqnarray*}
So, we can finish the second case by showing that
$$\||z|^{-1} {\bf I}_1^D f\|^r_{L^r(D)} 
= \int_D \frac{|z|^{-r} |z|^r}{2\pi} \left|\int_D \frac{f(\zeta)}{\zeta^s} \cdot \frac{d\zeta\wedge d\o{\zeta}}{\zeta^{1-s}(\zeta-z)} \right|^r dV(z)
\lesssim \|f\|^r_{|z|^s L^p(D)},$$
but that is an easy consequence of Theorem \ref{thm:young}.
In fact, consider the kernel
$$\Phi(\zeta,z) = \frac{1}{\zeta^{1-s}(\zeta-z)},$$
which appears in the situation that we are considering now. $\Phi$ satisfies the assumptions of Theorem \ref{thm:young}
with our choices of $r$, $p$ and $t$ because of \eqref{eq:ref4}, and so the case $\wt{k}(p,s)=k_1(p,s)=1$ and $s>0$ is settled.\\

Now, assume that $k(p,s)=k_1(p,s)=1$ and $s=0$, which implies that \linebreak $p=(1+s)p >2$.
This situation is well-known: Namely,
the Operator ${\bf I}_0^D$ maps continuously ${\bf I}_0^D: L^p(D)\rightarrow C^{\eta}(D)$ for all $2<p\leq\infty$, $0\leq \eta<1-2/p$.

For a proof, see \cite{He1}, Hilfssatz 15.
So, choose $\eta=1-2/p - \epsilon/2$. Then we deduce that there exists a constant $C(p,\eta)>0$
such that
$$|{\bf I}_1^D f(z)| = |{\bf I}_0^D f(z) - {\bf I}_0^D f(0)| \leq C(p,\eta) \|f\|_{L^p(D)} |z|^\eta \lesssim |z|^{1-2/p-\epsilon/2} \|f\|_{L^p(D)}.$$
But then (for $p<\infty$):
\begin{eqnarray*}
\|{\bf I}_1^D f\|^p_{|z|^{1-\epsilon} L^p(D)} &=& \int_D \left||z|^{\epsilon-1} {\bf I}_1^D f(z)\right|^p dV(z)\\
&\lesssim& \|f\|^p_{L^p(D)} \int_D |z|^{\frac{\epsilon}{2}p -2} dV(z) \lesssim \|f\|^p_{L^p(D)},
\end{eqnarray*}
completing the proof of Theorem \ref{thm:weighted} for $k(p,s)=k_1(p,s)=1$.
The case $p=\infty$ is also clear.\\

It remains to consider the situation $k(p,s)=k_1(p,s)=2$, thus
$sp > 2$,
and $s>0$ if $p=\infty$. Again, we just use the H{\"o}lder regularity of ${\bf I}_0^D$.
Let $q<2/(1-s)$. Then $|\zeta|^{s-1} \in L^q(D)$, and it follows that $\zeta^{-1} f\in L^r(D)$ for all
$1\leq r\leq\infty$ satisfying
$$\frac{1}{r} = \frac{1}{p} + \frac{1}{q} > \frac{1}{p} + \frac{1-s}{2} = \frac{2+(1-s)p}{2p},$$
because
$$\|\zeta^{-1} f\|_{L^r(D)} \leq \|\zeta^{-s} f\|_{L^p(D)} \|\zeta^{s-1}\|_{L^q(D)} = \|\zeta^{s-1}\|_{L^q(D)} \|f\|_{|z|^s L^p(D)}.$$
So, we can use that
\begin{eqnarray}\label{eq:new}
\|\zeta^{-1} f\|_{L^r(D)} \lesssim \|f\|_{|z|^s L^p(D)} \ \ \mbox{ for all }\ \ 1\leq r < \frac{2p}{2+(1-s)p}.
\end{eqnarray}
This implies that ${\bf I}^D_0(\zeta^{-1} f)=z^{-1}{\bf I}^D_1 f$ is in $C^\eta$ for all 
$$0\leq \eta < 1 - \frac{2 + (1-s)p}{p} = s-\frac{2}{p} < 1.$$
Choose $\eta=s-2/p-\epsilon/2$. So,
\begin{eqnarray*}
|z^{-1} {\bf I}^D_2 f(z)| &=& | {\bf I}_0^D(\zeta^{-1} f)(z) - {\bf I}_0^D (\zeta^{-1} f)(0)|\\
&\lesssim& |z-0|^\eta \|\zeta^{-1}f\|_{L^r(D)} \lesssim |z|^\eta \|f\|_{|z|^sL^p(D)}.
\end{eqnarray*}
But then (for $p<\infty$):
\begin{eqnarray*}
\|{\bf I}_2^D f\|^p_{|z|^{s+1-\epsilon} L^p(D)} &=& \int_D \left||z|^{\epsilon-1-s} {\bf I}_2^D f(z)\right|^p dV(z)\\
&\lesssim& \|f\|^p_{|z|^s L^p(D)} \int_D |z|^{p(\epsilon-s + \eta)} dV(z) \\
&=& \|f\|^p_{|z|^s L^p(D)} \int_D |z|^{\frac{\epsilon}{2}p -2} dV(z) \lesssim \|f\|^p_{|z|^s L^p(D)},
\end{eqnarray*}
completing the proof of Theorem \ref{thm:weighted} for $k(p,s)=k_1(p,s)=2$,
and again the case $p=\infty$ is clear, as well.

\section{The Basic $\dq$-Homotopy Formula on Product Domains}\label{sec:homotopy}

It is well known that the Grothendieck-Dolbeault Lemma can be proved by an inductional procedure
using only the inhomogeneous Cauchy Integral Formula in one complex variable (see e.g. \cite{GrRe}).
In this section, we will extend this to a kind of homotopy formula with error terms for the $\dq$-equation
on product domains that involves only Cauchy's Integral Formula.

\begin{defn}
We call $P=D_1\times \cdots \times D_n\subset \C^n$ a piecewise smooth product domain
if the $D_j$ are piecewise smooth domains 
in the complex plane for all $j=1, ..., n$.
Let
$$b_e P = \o{D_1}\times\cdots\times b D_e \times \cdots \times \o{D_n} = 
\{z=(z_1, ..., z_n)\in\o{P}: z_e\in bD_e\},$$
where $b D_e$ is the boundary of $D_e$.
\end{defn}

For the induction, we need the following

\begin{defn}
For $A\subset \C^n$ and $e\in \Z$, $1\leq e \leq n$, let $\Gamma_e(A)\subset \bigcup_{q=1}^n L^1_{(0,q),loc}(A)$ 
be the subset of locally integrable $(0,q)$-forms $\omega$ with $\dq \omega \in L^1_{(0,q+1),loc}(A)$ such that
$\omega$ and $\dq \omega$ do not contain $d\o{z_{e+1}}$, .... $d\o{z_{n}}$ .
For arbitrary $\eta \in L^1_{0,q}(A)$, we write
$$\eta=\eta_e + \eta'$$
for the unique representation where $\eta'$ does not contain $d\o{z_e}$, and $\eta_e$ contains $d\o{z_e}$ necessarily. 
\end{defn}

For the integration, we will use the following operators:

\begin{defn}
Let $P=D_1\times\cdots\times D_n\subset\subset\C^n$ be a piecewise smooth, bounded product domain.
Moreover, let $k\in\Z$, $e,q\in \Z$ with $1\leq e,q \leq n$, and let
$\omega$ be a $(0,q)$-form with measurable coefficients on $P$, given in multi-index notation as
$$\omega=\omega_e+\omega'=\sum_{\substack{|J|=q-1\\e\notin J}} a_{eJ} d\o{z_e}\wedge d\o{z_J} 
+ \sum_{\substack{|L|=q\\e\notin L}} a_L d\o{z_L}.$$
For $z\in \C^n$, let
$${\bf I}^{P_e}_k \omega (z):= \frac{z^k_e}{2\pi i} 
\sum_{\substack{|J|=q-1\\e\notin J}}  \int_{\zeta_e\in D_e} a_{eJ}(..., z_{e-1},\zeta_e,z_{e+1},...) \frac{d\zeta_e\wedge d\o{\zeta_e}}{\zeta_e^k(\zeta_e-z_e)} d\o{z_J},$$
provided the integral exists. 
\end{defn}

Note that ${\bf I}^{P_e}_0$
defines a bounded linear operator ${\bf I}^{P_e}_0: C^0_{0,q}(\o{P}) \rightarrow C^0_{0,q-1}(\o{P})$.\\

Moreover, we will have to describe some error terms which will be handled by integration over parts of the boundary.

For the representation of the error terms, we use:

\begin{defn}
Let $P=D_1\times\cdots\times D_n\subset\subset\C^n$ be a piecewise smooth, bounded product domain.
Moreover, let $k\in\Z$, $e,q\in \Z$ with $1\leq e,q \leq n$, and let
$\omega$ be a measurable $(0,q)$-form on $\o{P}$,
given in multi-index notation as
$$\omega=\omega_e+\omega'=\sum_{\substack{|J|=q-1\\e\notin J}} a_{eJ} d\o{z_e}\wedge d\o{z_J} 
+ \sum_{\substack{|L|=q\\e\notin L}} a_L d\o{z_L},$$
For $z\in \o{P}\setminus b_eP$, let
$${\bf R}^{P_e}_k \omega (z):= \frac{z_e^k}{2\pi i} \sum_{\substack{|L|=q\\e\notin L}}  
\int_{\zeta_e\in bD_e} a_{L}(..., z_{e-1},\zeta_e,z_{e+1},...) \frac{d\zeta_e}{\zeta_e^k(\zeta_e-z_e)} d\o{z_L},$$
provided the integral exists.
\end{defn}

Note that this Definition gives
a linear operator ${\bf R}^{P_e}_0: C^0_{0,q}(\o{P}) \rightarrow C^0_{0,q}(\o{P}\setminus b_e P)$,
for example. We will use the classical inhomogeneous Cauchy Integral Formula: 
Let $D\subset\subset \C$ be a bounded piecewise smooth domain, and assume that $f\in C^1(\o{D})$.
Then:
\begin{eqnarray}\label{eq:cif9}
f(z) = \frac{1}{2\pi i} \int_{D} \frac{\partial f}{\partial \o{\zeta}}(\zeta) \frac{d\zeta\wedge d\o{\zeta}}{\zeta-z} 
+ \frac{1}{2\pi i} \int_{bD} f(\zeta)\frac{d\zeta}{\zeta-z}
\end{eqnarray}
for all $z\in D$.
Now we are in the position to construct a homotopy formula with error terms which will have interesting applications.

\begin{lem}\label{lem:step1}
Let $P=D_1\times\cdots\times D_n\subset\subset\C^n$ be a piecewise smooth, bounded product domain,
$1\leq q \leq e\leq n$, and $\omega\in C^\infty_{0,q}(\o{P})$ a smooth $(0,q)$-form on a neighborhood of $\o{P}$.
Then:
\begin{eqnarray}\label{eq:homotopy}
\omega  = \dq {\bf I}_0^{P_e} \omega + {\bf I}_0^{P_e} \dq \omega + {\bf R}_0^{P_e} \omega\ \ \ \mbox{ on } \o{P}\setminus b_e P,
\end{eqnarray}
where ${\bf I}_0^{P_e}\omega\in C^\infty_{0,q-1}(\o{P}\setminus b_e P)$ 
and ${\bf I}_0^{P_e} \dq \omega, {\bf R}_0^{P_e} \omega \in C^\infty_{0,q}(\o{P}\setminus b_e P)$.
If moreover $\omega\in \Gamma_e(\o{P})$, then ${\bf R}_0^{P_e}\omega \in \Gamma_{e-1}(\o{P}\setminus b_e P)$,
which has the meaning that ${\bf R}_0^{P_q} \omega=0$ if $\omega\in \Gamma_q(\o{P})$.
If $\dq\omega=0$, then $\dq {\bf R}_0^{P_e} \omega=0$, too.
\end{lem}

\begin{proof}
Let us first consider the case when $\omega$ is of the form
$$\theta = a\ d\o{z_e}\wedge d\o{z_L} \in C^\infty_{0,q}(\o{P}).$$
So, assume that $\theta$ has the same properties as $\omega$.
Then:
\begin{eqnarray*}
\dq {\bf I}_0^{P_e} \theta &=& \theta + \frac{1}{2\pi i} \sum_{j\neq e} \int_{\zeta_e \in D_e} \frac{\partial a}{\partial\o{z_j}} 
\cdot\frac{d\zeta_e\wedge d\o{\zeta_e}}{\zeta_e - z_e} d\o{z_j}\wedge d\o{z_L}\\
&=& \theta - {\bf I}_0^{P_e} \dq \theta.
\end{eqnarray*}
Summing up the components of $\omega_e$ (which all look like $\theta$),
we conclude:
\begin{eqnarray}\label{eq:new1}
\omega_e &=& \dq {\bf I}_0^{P_e} \omega_e + {\bf I}_0^{P_e} \dq \omega_e = \dq {\bf I}_0^{P_e} \omega + {\bf I}_0^{P_e} \dq\omega_e
= \dq {\bf I}_0^{P_e} \omega + {\bf I}_0^{P_e} \dq\omega - {\bf I}_0^{P_e} \dq \omega'.
\end{eqnarray}
By use of the inhomogeneous Cauchy Integral Representation Formula \eqref{eq:cif9}
applied to the coefficients of $\omega'=\omega-\omega_e$ we also have:
$$\omega' = {\bf I}_0^{P_e} \dq \omega' + {\bf R}_0^{P_e} \omega'.$$
Inserting this to \eqref{eq:new1}, we arrive at
\begin{eqnarray*}
\omega_e &=& \dq {\bf I}_0^{P_e} \omega + {\bf I}_0^{P_e} \dq \omega - \omega' + {\bf R}_0^{P_e} \omega'\\
&=& \dq {\bf I}_0^{P_e} \omega + {\bf I}_0^{P_e} \dq \omega - \omega' + {\bf R}_0^{P_e} \omega,
\end{eqnarray*}
which proves \eqref{eq:homotopy}. Now, $\omega \in C^\infty_{0,q}(\o{P})$ implies that
${\bf I}_0^{P_e}\omega\in C^\infty_{0,q-1}(\o{P}\setminus b_e P)$ and 
${\bf I}_0^{P_e} \dq \omega, {\bf R}_0^{P_e} \omega \in C^\infty_{0,q}(\o{P}\setminus b_e P)$
by well-known properties of the Cauchy Integrals (see \cite{Ra}, Lemma IV.1.14).
If $\omega\in\Gamma_e(\o{P})$, then it is clear that ${\bf R}_0^{P_e} \omega\in\Gamma_e(\o{P}\setminus b_e P)$. 
But moreover, ${\bf R}_0^{P_e}\omega$ is holomorphic in $z_e$, and that implies that ${\bf R}_0^{P_e} \omega \in\Gamma_{e-1}(\o{P}\setminus b_e P)$.
The last statement is clear because $\dq\omega=0$ implies that 
$$\omega = \dq {\bf I}_0^{P_e} \omega + {\bf R}_0^{P_e} \omega.$$
\end{proof}

Applying Lemma \ref{lem:step1} to a $(0,q)$-form $n-q+1$ times yields a nice representation formula,
which we call a basic $\dq$-homotopy formula on product domains.
Let us define the inductively given integral operators:

\begin{defn}
Let $P=D_1\times\cdots\times D_n\subset\subset\C^n$ be a piece-wise smooth, bounded product domain, and
let $1\leq q \leq e\leq n$ and $\omega\in C^\infty_{0,q}(\o{P})$. Then we define
\begin{eqnarray*}
{\bf S}_q^P \omega &:=& \sum_{k=q}^n {\bf I}_0^{P_k}\ {\bf R}_0^{P_{k+1}}\ {\bf R}_0^{P_{k+2}} \cdots {\bf R}_0^{P_n}\ \omega,\\
{\bf T}_{q+1}^P \dq \omega &:=& {\bf I}_0^{P_n} \dq\omega + \sum_{k=q}^{n-1} {\bf I}_0^{P_k}\ \dq {\bf R}_0^{P_{k+1}}\ \dq {\bf R}_0^{P_{k+2}} \cdots 
\ \dq {\bf R}_0^{P_n}\ \omega\\
&=& {\bf I}_0^{P_n} \dq\omega + \sum_{k=q}^{n-1} {\bf I}_0^{P_k}\ \dq {\bf R}_0^{P_{k+1}}\ \dq {\bf R}_0^{P_{k+2}} \cdots 
\ \dq {\bf R}_0^{P_{n-1}} (\dq\omega - \dq {\bf I}_0^{P_n} \dq\omega).
\end{eqnarray*}
Note that ${\bf S}_q^P \omega \in C^\infty_{0,q-1}(P)$ and ${\bf T}_{q+1}^P \dq \omega \in C^\infty_{0,q} (P)$.
We have used 
$$\dq {\bf R}_0^{P_n} \omega = \dq\omega - \dq {\bf I}_0^{P_n} \dq\omega,$$
which follows from applying $\dq$ to \eqref{eq:homotopy}.
\end{defn}

Now we get the following $\dq$-homotopy formula:

\begin{thm}\label{thm:homotopy1}
Let $P=D_1\times\cdots\times D_n\subset\subset\C^n$ be a piece-wise smooth, bounded product domain, and
let $1\leq q \leq n$ and $\omega\in C^\infty_{0,q}(\o{P})$. Then:
\begin{eqnarray*}
\omega = \dq {\bf S}_q^P \omega + {\bf T}_{q+1}^P \dq \omega,
\end{eqnarray*}
where ${\bf S}_q^P \omega \in C^\infty_{0,q-1}(P)$ and ${\bf T}_{q+1}^P \dq \omega \in C^\infty_{0,q} (P)$. If $\dq\omega=0$, then $\omega=\dq {\bf S}_q^P \omega$.
\end{thm}

\begin{proof}
Applying Lemma \ref{lem:step1} to $\omega\in C^\infty_{0,q}(\o{P}) = C^\infty_{0,q}(\o{P})\cap \Gamma_n(\o{P})$ (trivially) yields
$$\omega = \dq {\bf I}_0^{P_n} \omega + {\bf I}_0^{P_n} \dq\omega + {\bf R}_0^{P_n} \omega,$$
where ${\bf R}_0^{P_n} \omega \in C^\infty_{0,q}(\o{P}\setminus b_n P)\cap \Gamma_{n-1}(\o{P}\setminus b_n P)$.
Now, using Lemma \ref{lem:step1} again, we have
$${\bf R}_0^{P_n} \omega = \dq {\bf I}_0^{P_{n-1}} {\bf R}_0^{P_n} \omega + {\bf I}_0^{P_{n-1}} \dq {\bf R}_0^{P_n} \omega + {\bf R}_0^{P_{n-1}} {\bf R}_0^{P_n} \omega,$$
where ${\bf R}_0^{P_{n-1}} {\bf R}_0^{P_n} \omega \in C^\infty_{0,q}(\o{P}\setminus (b_n P \cup b_{n-1} P))\cap \Gamma_{n-2}(\o{P}\setminus (b_n P\cup b_{n-1} P))$.
For this, the statement (not the proof) of Lemma \ref{lem:step1} has to be modified slightly, namely by replacing $\o{P}$ by $\o{P}\setminus b_n P$,
which does not really make a difference when integrating in the $z_{n-1}$-direction. Go on inductively. The induction ends at:
$${\bf R}_0^{P_{q+1}} \cdots {\bf R}_0^{P_n} \omega = \dq {\bf I}_0^{P_q}  {\bf R}_0^{P_{q+1}} \cdots {\bf R}_0^{P_n} \omega 
+ {\bf I}_0^{P_q} \dq {\bf R}_0^{P_{q+1}} \cdots \dq {\bf R}_0^{P_n} \omega,$$
because ${\bf R}_0^{P_{q+1}}\cdots {\bf R}_0^{P_n}\omega \in \Gamma_q$, and that implies ${\bf R}_0^{P_q} {\bf R}_0^{P_{q+1}}\cdots {\bf R}_0^{P_n} \omega =0$.
\end{proof}

\section{Local $L^p$-Solution of $\dq$ with Weights According to \\ Normal Crossings}\label{sec:solution}

We will now use a modification of Lemma \ref{lem:step1} in order to construct an $L^p$-solution operator
for the weighted $\dq_k$-equation similarly to the derivation of the homotopy formula Theorem \ref{thm:homotopy1}. 
One has to face the problem that Lemma \ref{lem:step1} can not be extended to $L^p$-forms directly
because there are boundary integrals involved. 
So, we simply restrict to the case where $\omega$ is compactly supported.
The resulting statement will be sufficient to deduce a local solution operator.
For the choice of the right weight factors in the occurring Cauchy Integrals,
we have to use the $\dq$-weights that we introduced in Definition \ref{defn:k} and Definition \ref{defn:wtk}.
Here now the basic ingredient for the construction of an $L^p$-solution operator:

\begin{lem}\label{lem:step2}
Let $P=D_1\times\cdots\times D_n\subset\subset\C^n$ be a 
bounded product domain, \linebreak
$1\leq q \leq e\leq n$, $1\leq p\leq \infty$, $s\in \R^n$, $\omega \in |z|^s L^p_{0,q}(P)$ 
with compact support in $D_e$,
and $k=k(p,s)$ the  $\dq$-weight of $(p,s)$ according to Definition \ref{defn:k}.
Let $m\in\Z^n$ with $m\leq k(p,s)$, 
and assume that
$\dq_{m} \omega \in |z|^s L^p_{0,q+1}(P)$.
Then:
\begin{eqnarray}\label{eq:homotopy2}
\omega  = \dq_{m} {\bf I}_{m_e}^{P_e} \omega + {\bf I}_{m_e}^{P_e} \dq_{m} \omega
\end{eqnarray}
in the sense of distributions.
\end{lem}

\begin{proof}
Let $r=k-m\in\Z^n$. Then $r \geq 0$.
By Definition \ref{defn:k} it follows that
\begin{eqnarray*}
z^{-m} \omega = z^r z^{-k}\omega\in L^1_{0,q}(P),
\end{eqnarray*}
and
\begin{eqnarray*}
\dq(z^{-m}\omega) = z^{-m} \dq_{m} \omega = z^r z^{-k} \dq_m \omega \in L^1_{0,q+1}(P).
\end{eqnarray*}
That shows 
that we are in a nice position to approximate $z^{-m} \omega$ by smooth forms.\\

It is well known (see e.g. \cite{LiMi}, Theorem V.2.6) that there exists a sequence of forms $\{g_j\}\subset C^\infty_{0,q}(P)$
such that
\begin{eqnarray*}
\lim_{j\rightarrow\infty} g_j &=& z^{-m} \omega \ \ \ \mbox{ in } \ L^1_{0,q}(P),\\
\lim_{j\rightarrow\infty} \dq g_j &=& \dq(z^{-m}\omega)=z^{-m}\dq_{m} \omega \ \ \ \mbox{ in } \ L^1_{0,q+1}(P).
\end{eqnarray*}
Because the sequence $\{g_j\}$ is constructed by convolution with a Dirac sequence, it is clear that we can assume
that the $g_j$ have compact support in $D_e$, as well. Lemma \ref{lem:step1} yields:
\begin{eqnarray*}
g_j = \dq {\bf I}^{P_e}_0 g_j + {\bf I}_0^{P_e} \dq g_j \ \ \ \mbox{ for all } j.
\end{eqnarray*}
This implies that
\begin{eqnarray}\label{eq:cif2}
z^{-m} \omega = \dq {\bf I}^{P_e}_0 (\zeta^{-m} \omega) + {\bf I}_0^{P_e} (\zeta^{-m} \dq_{m} \omega)
\end{eqnarray}
in the sense of distributions
by the regularity properties of the Cauchy-Integral,
where we use
$$\zeta=(z_1, ..., z_{e-1}, \zeta_e, z_{e+1}, ..., z_n).$$
But this implies that \eqref{eq:cif2} is equivalent to:
\begin{eqnarray*}
\omega &=& z^{m} \dq\big( z^{-m} {\bf I}^{P_e}_{m_e} (\omega)\big) + z^{m} z^{-m} {\bf I}^{P_e}_{m_e} (\dq_{m} \omega)\\
&=& \dq_{m} {\bf I}^{P_e}_{m_e} (\omega) +  {\bf I}^{P_e}_{m_e} (\dq_{m} \omega).
\end{eqnarray*}
\end{proof}

We will now use Lemma \ref{lem:step2} to construct a local solution operator for the $\dq_k$-equation.
One could go for a homotopy formula again, but the statement wouldn't be so nice and we actually don't need it.
So,
let $P=D_1\times\cdots \times D_n$ and \linebreak $Q=G_1\times \cdots \times G_n$ be two
bounded product domains in $\C^n$ such that 
$$Q\subset\subset P\subset\subset \C^n,$$
and choose smooth cut-off functions $\chi_j\in C^\infty_{cpt}(D_j)$ such that $\chi_j\equiv 1$
in a neighborhood of $G_j$ and $0\leq \chi_j \leq 1$. Let
$$S_j := D_1\times\cdots\times D_{j-1}\times \big(D_j\setminus G_j\big)\times D_{j+1}\times \cdots \times D_n.$$
Note that 
$$\supp\dq\chi_j \subset S_j$$
if we interpret $\chi_j$ as a function on $P$ (only depending on $z_j$).
Moreover, let $1\leq p \leq \infty$,
and $s=(s_1, ..., s_n)\in \R^n$ be a multi-index.
We will treat two slightly different situations in one.
So, let $c\in \Z^n$ be one of the following two weights:
$c(p,s) = k = k(p,s)$
the $\dq$-weight of $(p,s)$, or
$c(p,s) = \wt{k} = \wt{k}(p,s)$
the modified $\dq$-weight of $(p,s)$.
Let
$$s^+:=(s_1, ..., s_n+1-\epsilon)$$
for an arbitrary $\epsilon>0$ in the first case, or $s^+:=s$ in the second case.\\

Now, let $\omega\in |z|^s L^p_{0,q}(P)$ such that
$$\dq_{c}\omega=0$$
in the sense of distributions.
Our aim is to find a solution $\eta \in |z|^{s^+} L^p_{0,q-1}(P)$ such that
$$\dq_c \eta= \dq_{\wt{k}} \eta=\omega$$
in the sense of distributions on $Q$ (note that $\dq_c \eta=\dq_{\wt{k}} \eta$ because $\wt{k}\leq c$).
We will construct $\eta$ by induction. For this, set 
$$\omega^n := \omega.$$
Now, for $j=n, ..., q$ we define inductively
\begin{eqnarray*}
\eta^j &:=& {\bf I}^{P_j}_{c_j} \big(\chi_j \omega^j\big),\\
\omega^{j-1} &:=& {\bf I}^{P_j}_{c_j} \big(\dq \chi_j \wedge  \omega^j\big),\\
\theta^j &:=& {\bf I}^{P_j}_{c_j} \big( \chi_j \wedge \dq_c \omega^j\big).
\end{eqnarray*}
We will see later (Lemma \ref{lem:properties}) that this is in fact well defined. 
Note that
$$\omega^{j-1} + \theta^j = {\bf I}^{P_j}_{c_j} \big( \dq_c( \chi_j \omega^j)\big),$$
because
$$\dq_c (\chi_j\omega^j) = \dq \chi_j\wedge \omega^j + \chi_j\wedge \dq_c \omega^j.$$
Let
$$\eta:= \sum_{j=q}^n \eta^j.$$
We claim that $\eta$ has the desired properties. For the proof, we need another class of forms:

\begin{defn}
For $A\subset \C^n$, $e\in \Z$, $1\leq e \leq n$, let $\wt{\Gamma}_e(A)$ 
be the set of measurable $(0,q)$-forms on $A$
which do not contain $d\o{z_{e+1}}$, .... $d\o{z_{n}}$.
\end{defn}

Now we can collect properties of the forms occurring in the inductional procedure:

\begin{lem}\label{lem:properties}
For $j=n-1, ..., q$, it follows from the previous construction that
\begin{eqnarray}\label{eq:ind1}
\omega^j &\in& |z|^{s^+} L^p_{0,q}(P) \bigcap \wt{\Gamma}_j(P),\\
\dq_c \omega^{j}  &=& \dq\chi_{j+1} \wedge \omega^{j+1} - {\bf I}^{P_{j+1}}_{c_{j+1}} \big( \dq\chi_{j+1}\wedge\dq_c \omega^{j+1}\big)
\in |z|^{s} L^p_{0,q+1}(P),\label{eq:ind2}\\
\supp(\dq_c \omega^j) &\subset& \bigcup_{l=j+1}^n \supp \dq\chi_l \subset \bigcup_{l=j+1}^n S_l,\label{eq:ind3}\\
\supp \theta^j &\subset& \bigcup_{l=j+1}^n \supp \dq\chi_l \subset \bigcup_{l=j+1}^n S_l.\label{eq:ind4}
\end{eqnarray}
Together with $\omega^n=\omega$ and $\dq_c\omega^n=\dq_c\omega=0$,
this implies that $\eta^j$, $\omega^{j-1}$ and $\theta^j$ are well-defined for all $j=n, ..., q$.
\end{lem}

\begin{proof}
First of all, note that $\omega^n=\omega$ implies that
$$\omega^{n-1}={\bf I}^{P_n}_{c_n} (\dq\chi_n\wedge \omega^n) \in |z|^{s^+} L^p_{0,q}(P) \bigcap \wt{\Gamma}_{n-1}(P)$$
by the Definition of ${\bf I}^{P_n}_{c_n}$ and Theorem \ref{thm:weighted} or Corollary \ref{cor:weighted}.
Applying Lemma \ref{lem:step2} to $\dq\chi_n\wedge\omega^n$ leads to
\begin{eqnarray*}
\dq \chi_n\wedge \omega^n &=& \dq_c {\bf I}^{P_n}_{c_n} (\dq \chi_n\wedge \omega^n) + {\bf I}^{P_n}_{c_n} (\dq\chi_n\wedge \dq_c \omega^n)\\
&=& \dq_c \omega^{n-1} + {\bf I}^{P_n}_{c_n} (\dq\chi_n\wedge \dq_c \omega^n) = \dq_c \omega^{n-1},
\end{eqnarray*}
which show \eqref{eq:ind1} and \eqref{eq:ind2} in the case $j=n-1$ (recall $\dq_c \omega^n=\dq_c \omega=0$). But
$$\dq_c \omega^{n-1} = \dq \chi_n\wedge\omega^n$$
also implies \eqref{eq:ind3} and \eqref{eq:ind4} for $j=n-1$, because
$$\theta^{n-1} = {\bf I}^{P_{n-1}}_{c_{n-1}}\big(\chi_{n-1}\wedge \dq_c \omega^{n-1}\big).$$
Now, assume that the statement of the Lemma is true for an index $j\leq n-1$.
Then, the Definition of ${\bf I}^{P_j}_{c_j}$ and Theorem \ref{thm:weighted} (or Corollary \ref{cor:weighted}) show that
$$\omega^{j-1}={\bf I}^{P_j}_{c_j} (\dq\chi_j\wedge \omega^j) \in |z|^{s^+} L^p_{0,q}(P) \bigcap \wt{\Gamma}_{j-1}(P).$$
Lemma \ref{lem:step2} implies that
\begin{eqnarray*}
\dq \chi_j\wedge \omega^j &=& \dq_c {\bf I}^{P_j}_{c_j} (\dq \chi_j\wedge \omega^j) + {\bf I}^{P_j}_{c_j} (\dq\chi_j\wedge \dq_c \omega^j)\\
&=& \dq_c \omega^{j-1} + {\bf I}^{P_j}_{c_j} (\dq\chi_j\wedge \dq_c \omega^j),
\end{eqnarray*}
showing \eqref{eq:ind2}.
This in turn gives
\begin{eqnarray*}
\supp\dq_c \omega^{j-1} \subset \supp\dq\chi_j \bigcup \supp \dq_c \omega^j \subset \bigcup_{l=j}^n \supp\dq\chi_l,
\end{eqnarray*}
and the same is true for $\supp\theta^{j-1}$, because
$$\theta^{j-1} = {\bf I}^{P_{j-1}}_{c_{j-1}}\big(\chi_{j-1}\wedge \dq_c\omega^{j-1}\big).$$
\end{proof}

Let us recall that 
$$\eta^j =  {\bf I}^{P_j}_{c_j} \big(\chi_j \omega^j\big)$$
for $j=n, ..., q$. Now, it is a direct consequence of Theorem \ref{thm:weighted} (or Corollary \ref{cor:weighted}),
Lemma \ref{lem:properties} and the Definition of ${\bf I}^{P_j}_{c_j}$ that
$$\eta^j \in |z|^{s^+} L^p_{0,q-1}(P) \bigcap \wt{\Gamma}_{j-1}(P)$$
for all $j=n, ..., q$. Hence:
$$\eta=\sum_{j=q}^n \eta^j \in |z|^{s^+} L^p_{0,q-1}(P) \bigcap \wt{\Gamma}_{n-1}(P).$$

Moreover, Lemma \ref{lem:step2} implies that
\begin{eqnarray*}
\dq_c \eta^j &=& \dq_c {\bf I}^{P_j}_{c_j} \big(\chi_j \omega^j\big)
= \chi_j\omega^j - {\bf I}^{P_j}_{c_j} \big(\dq_c(\chi_j \omega^j)\big)\\
&=& \chi_j\omega^j -\omega^{j-1} - \theta^j
\end{eqnarray*}
for $j=n, ..., q+1$. 
In the case $j=q$, note that $\omega^q\in \wt{\Gamma}_q(P)$
contains $d\o{z_q}$ necessarily (so $\dq\chi_q\wedge\omega_q=0$),
such that
$$\dq_c(\chi_q \omega^q) = \chi_q\wedge \dq_c\omega^q,$$
giving
$\dq_c \eta^q = \chi_q\omega^q - \theta^q$.
So, let's have a look at
\begin{eqnarray*}
\dq_c \eta &=& \sum_{j=q}^n \dq_c \eta^j
= \sum_{j=q+1}^n \left( \chi_j\omega^j -\omega^{j-1} - \theta^j\right) + \chi_q\omega^q -\theta^q\\
&=& \chi_n\omega + \sum_{j=q}^{n-1} (\chi_j-1)\omega^j - \sum_{j=q}^n \theta^j.
\end{eqnarray*}
Note that $(\chi_j-1)|_Q\equiv 0$ by Definition and that $\theta^j|_Q\equiv 0$ 
by property \eqref{eq:ind4} in Lemma \ref{lem:properties}.
But that yields:
$$(\dq_c \eta)|_Q = (\chi_n\omega)|_Q = \omega.$$
It is clear from our construction that $\eta$ depends linearly on $\omega$,
and that this linear application maps continuously from $|z|^s L^p_{0,q}(P)$
into $|z|^{s^+} L^p_{0,q-1}(P)$. Writing down the explicit formula for that operator seems to be a little
messy, so, summing up, we conclude our main result:

\begin{thm}\label{thm:lp-solution}
Let $P=D_1\times\cdots \times D_n$ and $Q=G_1\times \cdots \times G_n$ be two
bounded product domains in $\C^n$ such that 
$$Q\subset\subset P\subset\subset \C^n.$$
Moreover, let $1\leq q\leq n$, $1\leq p \leq \infty$,
and $s=(s_1, ..., s_n)\in \R^n$ be a multi-index 
with 
$$c(p,s) = k(p,s)=(k_1(p,s_1), ..., k_n(p,s_n))$$
the $\dq$-weight of $(p,s)$  according to Definition \ref{defn:k}, or
$$c(p,s) = \wt{k}(p,s) = (\wt{k_1}(p,s_1), ..., \wt{k_n}(p,s_n)),$$
the modified $\dq$-weight of $(p,s)$ according to Definition \ref{defn:wtk}.
Let $s^+:=s$ in the second case, or $s^+:=(s_1, ..., s_{n-1}, s_n+1-\epsilon)$
for an arbitrary $\epsilon>0$ in the first case.
Then there exists a linear mapping ${\bf J}_{p,s^+,q}^{Q,P}$
which defines a bounded linear operator
$${\bf J}_{p,s^+,q}^{Q,P}: |z|^s L^p_{0,q} (P) \cap \ker \dq_c  \rightarrow |z|^{s^+} L^p_{0,q-1}(P)\ ,$$
such that
$$\dq_c {\bf J}_{p,s^+,q}^{Q,P} \omega = \dq_{\wt{k}} {\bf J}_{p,s^+,q}^{Q,P} \omega =  \omega$$
in the sense of distributions on $Q$ ($\dq_c {\bf J}_{p,s^+,q}^{Q,P} \omega=\dq_{\wt{k}} {\bf J}_{p,s^+,q}^{Q,P} \omega$ because $\wt{k}\leq c$).
\end{thm}


\section{About Theorem \ref{thm:main}}

Theorem \ref{thm:main} is an easy consequence of the $\dq_k$-solution presented in Theorem \ref{thm:lp-solution}:
Let $1\leq p\leq\infty$, $s\in\R^n$, and $k(p,s)$ the $\dq$-weight of $(p,s)$ according to Definition \ref{defn:k}.
We have seen that
$$\mathcal{I}^k \OO \subset |z|^s \mathcal{L}^p_{0,0}$$
as a consequence of Lemma \ref{lem:kk}. Applying Theorem \ref{thm:lp-solution} with
$$c(p,s) = k(p,s)$$
shows that the sequence \eqref{eq:complex} is exact at $|z|^s\mathcal{L}^p_{0,q}$ for all $q\geq 1$
because clearly $|z|^{s^+}\mathcal{L}^p_{0,q} \subset |z|^s\mathcal{L}^p_{0,q}$.
It remains to show that it is exact at $|z|^s \mathcal{L}^p_{0,0}$.
So, let $f \in (|z|^s \mathcal{L}^p_{0,0})_w$ be a germ such that $\dq_k f=0$.
But by Definition \ref{defn:dqk}, $\dq_k f=0$ implies that 
$$z^{-k}f\in (\mathcal{L}^1_{0,0})_w\ \ \ \mbox{ and }\ \ \dq(z^{-k} f)=0.$$
So, it follows that $f \in (\mathcal{I}^k \OO)_w$.
Now, let $f\in (|z|^s \mathcal{L}^p_{0,q})_w$ and $\phi\in (\mathcal{C}^\infty)_w$. Then:
$$\phi\cdot f \in (|z|^s L^p_{0,q})_w,$$
and 
$$\dq_k( \phi\cdot f) = \phi \dq_k f + \dq \phi \wedge f \in (|z|^s L^p_{0,q+1})_w.$$
That shows that the resolution is fine and completes the proof of Theorem \ref{thm:main}.\\

It is worth mentioning that we have also provided all the necessary knowledge to study
the $\dq_{\wt{k}}$-equation on $|z|^s L^p$-forms, where $\wt{k}(p,s)$ is the modified $\dq$-weight.
Precisely:

\begin{thm}\label{thm:main'}
For $1\leq p\leq \infty$ and $s\in \R^n$, let $\wt{k}=\wt{k}(p,s)\in \Z^n$ be the modified $\dq$-weight of $(p,s)$
according to Definition \ref{defn:wtk}.
Then:
\begin{eqnarray}\label{eq:complex'}
0 \rightarrow \mathcal{I}^{\wt{k}}\OO \hookrightarrow
 |z|^s \mathcal{L}^p_{0,0} \xrightarrow{\ \dq_{\wt{k}}\ }
 |z|^s \mathcal{L}^p_{0,1} \xrightarrow{\ \dq_{\wt{k}}\ }
\cdots \xrightarrow{\ \dq_{\wt{k}}\ }
 |z|^s \mathcal{L}^p_{0,n} \rightarrow 0
\end{eqnarray}
is an exact (and fine) resolution of $\mathcal{I}^{\wt{k}}\OO$.
\end{thm}

\begin{proof}
The proof is completely analogous to the proof of Theorem \ref{thm:main}.
Note that
$$\mathcal{I}^{\wt{k}}\OO \subset |z|^s \mathcal{L}^p_{0,0}$$
by the Definition \ref{defn:wtk} of the modified $\dq$-weight.
Theorem \ref{thm:lp-solution} with
$$c(p,s) = \wt{k}(p,s)$$
shows that the sequence \eqref{eq:complex'} is exact at $|z|^s\mathcal{L}^p_{0,q}$ for all $q\geq 1$,
and the rest of the proof goes through as above with $\wt{k}=\wt{k}(p,s)$ instead of $k=k(p,s)$.
\end{proof}


{\bf Acknowledgments}

\vspace{2mm}
This work was done while the author was visiting the University of Michigan at Ann Arbor,
supported by a fellowship within the Postdoc-Programme of the German Academic Exchange Service (DAAD).
The author would like to thank the SCV group at the University of Michigan for its hospitality.



\vspace{5mm}
\begin{thebibliography}{99999}













\bibitem[AHL]{AHL} {\sc J.\ M.\ Aroca, H.\ Hironaka, J.\ L.\ Vicente},
Desingularization theorems, {\em Mem. Math. Inst. Jorge Juan}, {\bf No. 30}, Madrid, 1977.


\bibitem[BiMi]{BiMi} {\sc E.\ Bierstone, P.\ Milman}, Canonical desingularization in characteristic zero by blowing-up
the maximum strata of a local invariant, {\em Inventiones Math.} {\bf 128} (1997), {\em no. 2}, 207--302.





\bibitem[CGM]{CGM} {\sc J.\ Cheeger, M.\ Goresky, R.\ MacPherson}, $L^2$-cohomology and intersection homology of singular algebraic varieties,
{\em Ann. Math. Stud.} {\bf 102} (1982), 303--340.







\bibitem[DFV]{DFV} {\sc K.\ Diederich, J.\ E.\ Forn{\ae}ss, S.\ Vassiliadou},
Local $L^2$ results for $\dq$ on a singular surface, {\em Math. Scand.} {\bf 92} (2003), 269--294.






\bibitem[Fo]{Fo} {\sc J.\ E.\ Forn{\ae}ss}, $L^2$ results for $\dq$ in a conic,
in {\em International Symposium, Complex Analysis and Related Topics, Cuernavaca, Operator Theory:
Advances and Applications} (Birkhauser, 1999).




\bibitem[FOV1]{FOV1} {\sc J.\ E.\ Forn{\ae}ss, N.\ {\O}vrelid, S.\ Vassiliadou},
Semiglobal results for $\dq$ on a complex space with arbitrary singularities,
{\em Proc. Am. Math. Soc.} {\bf 133} (2005), {\em no. 8}, 2377--2386.


\bibitem[FOV2]{FOV2} {\sc J.\ E.\ Forn{\ae}ss, N.\ {\O}vrelid, S.\ Vassiliadou},
Local $L^2$ results for $\dq$: the isolated singularities case,
{\em Internat. J. Math.} {\bf 16} (2005), {\em no. 4}, 387--418.








\bibitem[GrRe]{GrRe}{\sc H.\ Grauert, R.\ Remmert}, {\em Theory of Stein Spaces},
Springer-Verlag, Berlin, 1979.






\bibitem[Ha]{Ha} {\sc H.\ Hauser}, The Hironaka theorem on resolution of singularities,
{\em Bull. (New Series) Amer. Math. Soc.} {\bf 40}, no. 3, (2003), 323--403.


\bibitem[He1]{He1} {\sc T.\ Hefer}, Optimale Regularit\"atss\"atze f\"ur die $\dq$-Gleichung
auf gewissen konkaven und konvexen Modellgebieten,
{\em Bonner Math. Schr.} {\bf 301} (1997).


































\bibitem[LiMi]{LiMi} {\sc I.\ Lieb, J.\ Michel}, {\em The Cauchy-Riemann Complex,
Integral Formulae and Neumann Problem}, Vieweg, Braunschweig/Wiesbaden, 2002.







\bibitem[Oh]{Oh} {\sc T.\ Ohsawa}, Cheeger-Goresky-MacPherson's conjecture for the varieties with isolated singularities,
{\em Math. Z.} {\bf 206} (1991), 219--224.







\bibitem[PaSt1]{PaSt1} {\sc W.\ Pardon, M.\ Stern}, $L^2$-$\dq$-cohomology of complex projective varieties,
{\em J. Amer. Math. Soc.} {\bf 3} (1991), no. 3, 603--621.



\bibitem[PaSt2]{PaSt2} {\sc W.\ Pardon, M.\ Stern}, Pure Hodge structure on the $L^2$-cohomology of varieties with isolated singularities,
{\em J. reine angew. Math.} {\bf 533} (2001), 55--80.






\bibitem[Ra]{Ra} {\sc R.\ M.\ Range}, {\em Holomorphic Functions and Integral 
Representations in Several Complex Variables}, (Graduate Texts in 
Mathematics, Bd. 108), Springer, New York, 1986.










\bibitem[Ru2]{Rp2} {\sc J.\ Ruppenthal}, Zur Regularit\"at der Cauchy-Riemannschen Differentialgleichungen
auf komplexen R\"aumen, {\em Bonner Math. Schr.} {\bf 380} (2006).


\bibitem[Ru4]{Rp4}{\sc J.\ Ruppenthal}, About the $\dq$-equation at isolated singularities with regular exceptional set,
preprint 2007, {\sf arXiv:0803.0152}, to appear in {\em Int. J. Math.}.


\bibitem[Ru7]{Rp7}{\sc J.\ Ruppenthal}, The $\dq$-equation on homogeneous varieties with an isolated singularity,
preprint 2008, {\sf arXiv:0803.1188}, submitted.



\bibitem[RuZe]{RuZe} {\sc J.\ Ruppenthal, E.\ S.\ Zeron}, An explicit $\dq$-integration formula
for weighted homogeneous varieties, preprint 2008, {\sf arXiv:0803.0136}, to appear in {\em Mich. Math. J.}












\end{thebibliography}
\end{document}